\documentclass[reqno]{amsart}
\usepackage[utf8]{inputenc}
\usepackage{amsfonts, amsmath, amssymb, amsthm, graphicx, enumitem}
\usepackage[colorlinks=true,urlcolor=blue, citecolor=red,linkcolor=blue]{hyperref}
\usepackage[left=2cm, right=2cm, top=3cm]{geometry}
\newtheorem{theorem}{Theorem}

\newtheorem{definition}[theorem]{Definition}

\newtheorem{lemma}[theorem]{Lemma}
\newtheorem{proposition}[theorem]{Proposition}
\newtheorem{remark}[theorem]{Remark}
\usepackage{amsfonts}
\usepackage{dsfont}

\renewcommand{\d}{\displaystyle}
\usepackage{graphicx}
\usepackage{apptools}
\AtAppendix{\counterwithin{lem}{section}}
\newtheorem{lem}{Lemma}
\newcommand{\pts}[1]{\left(#1\right)}  
\newcommand{\cts}[1]{\left[#1\right]}                                  	  %%
\newcommand{\lvs}[1]{\left\{#1\right\}}                                	  %%
\newcommand{\abs}[1]{\left|#1\right|}                                  	  %%

\newcommand{\C}{\mathbb{C}} 
\newcommand{\R}{\mathbb{R}}
\newcommand{\A}{\mathcal{A}}
\newcommand{\al}{\alpha}

\newcommand{\U}{\mathcal{U}}
\newcommand{\OO}{\mathcal{O}}
\newcommand{\Fi}{\Phi}

\newcommand{\HH}{\mathcal{H}}

\begin{document}
	\title[Boundary controllability for a 1D degenerate parabolic equation]{Boundary controllability for a 1D degenerate parabolic equation with a Robin boundary condition}
	
	\author[L. Galo-Mendoza \and M. L\'opez-Garc\'ia]{Leandro Galo-Mendoza \and Marcos L\'opez-Garc\'ia}
	\address{Instituto de Matem\'{a}ticas-Unidad Cuernavaca\\
		Universidad Nacional Aut\'onoma de M\'exico\\
		Av. Universidad S/N\\
		Cuernavaca, Morelos, C.P. 62210\\
		M\'{e}xico.}
	\email{jesus.galo@im.unam.mx, marcos.lopez@im.unam.mx}
	
	
	\maketitle
	
	\begin{abstract}In this paper we prove the null controllability of a one-dimensional degenerate parabolic equation with a weighted Robin boundary condition at the left endpoint, where the potential has a singularity. We use some results from the singular Sturm-Liouville theory to show the well-posedness of our system. We obtain a spectral decomposition of a degenerate parabolic operator with Robin conditions at the endpoints, we use Fourier-Dini expansions and the moment method introduced by Fattorini and Russell to prove the null controllability and to obtain an upper estimate of the cost of controllability. We also get a lower estimate of the cost of controllability by using a representation theorem for analytic functions of exponential type.
	\end{abstract}

	\section{Introduction and main results}
	Let $T > 0$ and set $Q := (0, 1) \times (0, T )$. For $\alpha, \beta\in \R$ with $0\leq \al < 2$, consider the equation
	\begin{equation}\label{degeqn}
		u_t-(x^\al u_x)_x-\beta x^{\al -1}u_x-\frac{\mu}{x^{2-\al}} u =0  \text { in }Q,
	\end{equation}
	provided that $\mu\in\R$ satisfies
	\begin{equation}\label{mucon}
		-\infty<\mu<\mu(\al+\beta),\quad \text{where}\quad\mu(\delta):=\frac{(1-\delta)^2}{4}, \quad \delta\in\R.%\quad \text{and}
	\end{equation}
	In this work, we consider a weighted Robin boundary condition at the left endpoint of the form
	\[ \lim_{x\rightarrow 0^+}\pts{ax^{(\al+\beta-1)/2+\sqrt{\mu(\al+\beta)-\mu}}u(x,t)+x^{(\al+\beta+1)/2+\sqrt{\mu(\al+\beta)-\mu}}u_{x}(x,t)} =f(t),\]
	and a usual Robin boundary condition at the right endpoint of the form
	$$
	au(1,t)+u_{x}(1,t) =g(t),
	$$ 
	where
	\begin{equation}\label{a_const}
		a:=a(\al,\beta,\mu)=\dfrac{\al+\beta-1}{2}-\sqrt{\mu(\alpha+\beta)-\mu}.
	\end{equation} 
	
	The goal of this work is to prove the null controllability of the following system, with a control $f(t)\in L^2(0,T)$ acting  at the left endpoint,
	\begin{equation}\label{problem1Left}
		\left\{\begin{aligned}
			u_t-(x^\al u_x)_x-\beta x^{\al -1}u_x-\frac{\mu}{x^{2-\al}} u&=0 & & \text { in }Q, \\
			\cts{u(\cdot,t),x^{-a}}(0) = f(t),\quad au(1, t)+u_{x}(1,t)&=0 & & \text { on }(0, T), \\
			u(x, 0) &=u_{0}(x) & & \text { in }(0, 1),
		\end{aligned}\right.
	\end{equation}
	where our Lagrange form $[\cdot,\cdot ]$ is given by
	\begin{equation*}
		\cts{u,v}(x)=(upv'-vpu')(x), \quad \text{with  } p(x)=x^{\alpha+\beta}, \text{ and } ^\prime=\frac{d}{dx}.
	\end{equation*}
	Consider the weighted Lebesgue space $L^2_\beta(0,1):=L^2((0,1);x^\beta dx)$, $\beta\in\R$, endowed with the inner product
	\[\langle f,g\rangle_\beta:=\int_0^1 f(x)g(x)x^\beta dx,\]
	and its corresponding norm is denoted by $\|\cdot\|_{\beta}$.\\
	
	Here, we use some results from the singular Sturm-Liouville theory to see the well-posedness of the system (\ref{problem1Left}) with initial data in $L^2_\beta(0,1)$, although the solution $u(t)$ lives in an interpolation space $\mathcal{H}^{-s}$. We say the system (\ref{problem1Left}) is null controllable in $L^2_\beta(0,1)$ at time $T>0$ with controls in $L^2(0,T)$, if for any $u_0\in L^2_\beta(0,1)$ there exists $f\in L^2(0,T)$ such that the corresponding solution satisfies $u(\cdot,T)\equiv 0$.\\ 
	
	We are also interested in the behavior of the cost of the controllability. Consider the set of admissible controls given by
	\[U(T,\al,\beta,\mu,u_0):=\{f\in L^2(0,T): u\text{ is solution of the system (\ref{problem1Left}) that satisfies }u(\cdot,T)\equiv 0\}.\]
	
	If $X$ is a subspace in $L^2_\beta(0,1)$, we define the cost of controllability for initial data in $X$ as follows
	\[\mathcal{K}_X(T,\al,\beta,\mu):=\sup_{u_0\in X,\|u_0\|_\beta=1}\inf\{|f|_{L^2(0,T)}:f\in U(T,\al,\beta,\mu,u_0)\}.\]
	The main result of this work is the following.
	\begin{theorem}\label{Teo1}
		Let $T>0$, $0\leq\al<2$, $\beta\in \R$, and $\mu$ satisfying (\ref{mucon}). The next statements hold.
		\begin{enumerate}
			\item \textbf{Existence of a control} For any $u_0\in L^2_\beta(0,1)$ there exists a control $f \in L^2(0, T )$ such that the solution $u$ to (\ref{problem1Left}) satisfies $u(\cdot,T ) \equiv 0$.
			\item \textbf{Upper bound of the cost} There exists a constant $c>0$ such that for every $\delta\in (0,1)$ we have
			\[\mathcal{K}_{\Phi_0^\perp}(T,\alpha,\beta,\mu)\leq \frac{c M(T,\alpha,\nu,\delta)T^{1/2}}{(\nu+1)\kappa_\al^{5/2}}
			\exp\pts{-\frac{T}{2}\kappa_\alpha^2 j_{\nu+1,1}^2}.\]
			where 
			\begin{equation}\label{Nu}
				\kappa_\al:=\frac{2-\al}{2},\quad \nu=\nu(\al,\beta,\mu):=\sqrt{\mu(\al+\beta)-\mu}/\kappa_\al,\quad \Phi_0(x)=\sqrt{2(\nu+1)\kappa_\al}\,x^{-a},
			\end{equation}
			$j_{\nu+1,1}$ is the first positive zero of the Bessel function $J_{\nu+1}$ (defined in the Appendix), and
			\[M(T,\alpha,\nu,\delta)=\pts{1+\frac{1}{(1-\delta)\kappa_\alpha^2 T}}\cts{\exp\pts{\frac{1}{\sqrt{2}\kappa_\alpha}}+\frac{1}{\delta^5}\exp\pts{\frac{3}{(1-\delta)\kappa_\alpha^2 T}}}\exp\pts{-\frac{(1-\delta)^{3/2}T^{3/2}}{8(1+T)^{1/2}}\kappa_\alpha^3 j_{\nu+1,1}^2}.\]			
			\item \textbf{Lower bound of the cost} There exists a constant $c>0$ such that
			\[c\pts{1+\frac{j_{\nu+1,2}^2}{j_{\nu+1,1}^2}}\frac{2^{\nu}|J_{\nu}(j_{\nu+1,1})|\exp{\left(\pts{\frac{1}{2}-\frac{\log 2}{\pi}}j_{\nu +1,2}\right)}}{\Gamma(\nu+1)^{-1}\pts{2T\kappa_\al}^{1/2}(j_{\nu+1,1})^\nu}\exp\pts{-\pts{j_{\nu+1,1}^2+\frac{j_{\nu+1,2}^2}{2}}\kappa_\alpha^2T  }\leq \mathcal{K}_{L^2_\beta}(T,\alpha,\beta,\mu),\]
			where $j_{\nu+1,2}$ is the second positive zero of the Bessel function $J_{\nu+1}$.
		\end{enumerate}
		%where $j_{\nu,2}$ is the second positive zero of the Bessel function $J_\nu$.
	\end{theorem}				
	We also analyze the null controllability of a similar system but the control acting at the right endpoint,
	\begin{equation}\label{problem1Right}
		\left\{\begin{aligned}
			u_t-(x^\al u_x)_x-\beta x^{\al -1}u_x-\frac{\mu}{x^{2-\al}} u&=0 & & \text { in }Q, \\
			\cts{u(\cdot,t),x^{-a}}(0) = 0,\,\,\, au(1, t)+u_{x}(1,t)&= f(t) & & \text { on }(0, T), \\
			u(x, 0) &=u_{0}(x) & & \text { in }(0, 1).
		\end{aligned}\right.
	\end{equation}
	Consider the corresponding set of admissible controls
	\[\widetilde{U}(T,\al,\beta,\mu,u_0)=\{f\in L^2(0,T): u \text{ is solution of the system (\ref{problem1Right}) that satisfies } u(\cdot,T)\equiv 0\},\]
	and the cost of the controllability given by
	\[\widetilde{\mathcal{K}}_{X}(T,\al,\beta,\mu):=\sup_{u_0\in X,\|u_0\|_\beta=1}\inf\{\|f\|_{L^2(0,T)}:f\in \widetilde{U}(T,\al,\beta,\mu,u_0)\},\]
	where $X$ is a subspace in $L^2_{\beta}(0,1)$.
	\begin{theorem}\label{Teo2}
		Let $T>0$, $\beta \in \R$, $0\leq\al<2$, and $\mu$ satisfying (\ref{mucon}). The next statements hold.
		\begin{enumerate}
			\item \textbf{Existence of a control} For any $u_0\in L^2_\beta(0,1)$ there exists a control $f \in L^2(0, T )$ such that the solution $u$ to (\ref{problem1Right}) satisfies $u(\cdot,T ) \equiv 0$.
			\item \textbf{Upper bound of the cost} There exists a constant $c>0$ such that for every $\delta\in (0,1)$ we have
			\[\widetilde{\mathcal{K}}_{\Phi_0^\perp}(T,\alpha,\beta,\mu)\leq \frac{c M(T,\alpha,\nu,\delta)T^{1/2}}{\kappa_{\al}^{\nu+1}\Gamma(v+2)}
			\pts{\dfrac{2\nu+1}{4T\mathrm{e}}}^{(2\nu+1)/4} \exp\pts{-\frac{T}{4}\kappa_\alpha^2 j_{\nu+1,1}^2}.\]
			\item \textbf{Lower bound of the cost} There exists a constant $c>0$ such that
			\[c\pts{1+\frac{j_{\nu+1,2}^2}{j_{\nu+1,1}^2}}\frac{\exp{\left(\pts{\frac{1}{2}-\frac{\log 2}{\pi}}j_{\nu+1,2}\right)}}{\pts{2T\kappa_\alpha}^{1/2}}\exp\pts{-\pts{j_{\nu+1,1}^2+\frac{j_{\nu,2}^2}{2}}\kappa_\alpha^2 T}\leq \widetilde{\mathcal{K}}_{L^2_\beta}(T,\alpha,\beta,\mu).\]
		\end{enumerate}
	\end{theorem}
	%%%%%%%%%%%%%%%%%%%%%%%%%%%%%%%%%%%%%%%%%%%%%
	%%%%%%%%%%%%%%%%%%%%%%%%%%%%%%%%%%%%%%%%%%%%%%
	In \cite{GaloLopez} the authors prove the null controllability of the equation (\ref{degeqn}) with a weighted Dirichlet boundary condition at the left endpoint, provided that $\alpha+\beta<1$. In the case  $\alpha+\beta>1$, in \cite{GaloLopez2} they get the null controllability of the equation (\ref{degeqn}) with a weighted Neumann boundary condition at the left endpoint. They consider initial data in $L^2_\beta(0,1)$ in both cases. In these works the authors prove suitable versions of a Hardy inequality to assure the well-posedness of their systems, but in the case $\alpha+\beta=1$ is necessary to consider some results from the singular Sturm-Liouville theory, see \cite{GaloLopez2}. Here, we use that approach to show the well-posedness of our system.\\
	
	%%%%%%%%%%%%%%%%%%%%%%%%%%%%%%%%%%%%%%%%%%
	%%%%%%%%%%%%%%%%%%%%%%%%%%%%%%%%%%%%%%%%%
	
	This paper is organized as follows. Section 2 uses some results from the singular Sturm-Liouville to show that the operator $\mathcal{A}$ given in (\ref{A_ope}) is self-adjoint. There, we also use Fourier-Dini expansions to show that $\mathcal{A}$ is diagonalizable, this allows us to consider initial data in some interpolation spaces. Next, we introduce a notion of a weak solution for both systems and then show the well-posedness of these systems.\\
	
	In Section 3 we prove Theorem \ref{Teo1} by using the moment method introduced by Fattorini \& Russell. Here, the idea is to construct a biorthogonal sequence to a family of exponentials involving the eigenvalues of $\mathcal{A}$. To do this we use some results from complex analysis to construct a suitable complex multiplier. As a consequence, we get an upper estimate of the cost of the controllability. Finally, we use a representation theorem, Theorem \ref{repreTheo}, to obtain a lower estimate of the cost of the controllability.\\
	
	In Section 4 we proceed as before to solve the case when the control acts at the right endpoint.
	
	%%%%%%%%%%%%%%%%%%%%%%%%%%%%%%%%%%%%%%%%%%%%%%%%%%%%%%%%%%%%%%%%%%%%%%%%
	\section{Functional setting and well-posedness}
	Consider the differential expression $M$ defined by
	\[Mu=-(pu_x)_x+qu\]
	where $\d p(x) = x^{\al+\beta}, q(x) = -\mu x^{-2+\al+\beta}, w(x) = x^{\beta}$.\\
	
	Clearly,
	\begin{equation*}\label{Acond2}
		1/p, q, w \in L_{\text{loc}}(0,1),\quad p,w >0\text{ on } (0,1),
	\end{equation*} 
	thus $Mu$ is defined a.e. for functions $u$ such that $u, pu_x\in AC_{\text{loc}}(0,1)$, where $AC_{\text{loc}}(0,1)$ is the space of all locally absolutely continuous functions in $(0,1)$.\\
	
	Now we introduce the operator $\mathcal{A}$ given by 
	\begin{equation}\label{A_ope}
		\A u:=w^{-1}Mu=-(x^\al u_x)_x-\beta x^{\al -1}u_x-\frac{\mu}{x^{2-\al}} u .
	\end{equation}
	From the theory developed in \cite{Zettl} we can build a self-adjoint domain $D(\A)$ for the operator $\mathcal{A}$.\\% densely in $L^2_{\beta}(0,1)$ to the operator $\A$ defined by the problem in \ref{A_ope}.\\
	
	For $\mu$ satisfying (\ref{mucon}), $0\leq \al<2$, and $\beta\in\R$, we set
	\begin{equation*}\label{Dmax}
		D_{\max}:=\left\{u\in AC_{\text{loc}}(0,1)\, |\,pu_x\in AC_{\text{loc}}(0,1),\, u, \A u\in L^2_{\beta}(0,1)\right\},\quad\text{and}
	\end{equation*}
	\[D(\A):=\left\{\begin{aligned}
		\{u\in D_{\max}\,| \lim_{x\rightarrow 0^+}x^{(\al+\beta-1)/2+\sqrt{\mu(\al+\beta)-\mu}}u(x)=(au+u_{x})(1)=0\} & &\text{if } \sqrt{\mu(\alpha+\beta)-\mu}<\kappa_\al, \\
		\{u\in D_{\max}\,| (au+u_{x})(1)=0\}& & \text{if } \sqrt{\mu(\alpha+\beta)-\mu}\geq\kappa_\al.
	\end{aligned}\right.
	\]
	Recall that the Lagrange form associated with $M$ is defined as follows,
	\[[u,v]:=upv_x-vpu_x,\quad\text{for all}\quad u,v\in D_{\max}.\]
	
	The next result shows that $\mathcal{A}$ is a diagonalizable operator whose Hilbert basis of eigenfunctions can be written in terms of the function $x^{1/2+\nu}$, the Bessel function of the first kind $J_{\nu}$ and the corresponding positive zeros $j_{\nu+1,k}$, $k\geq 1$, of the Bessel function $J_{\nu+1}$, see the proof of Proposition \ref{basis}. In the appendix, we give some properties of Bessel functions of the first kind and their zeros.
	\begin{proposition}
		Let $0\leq \al<2$, $\beta\in\R$, $\mu < \mu(\alpha+\beta)$, and $\kappa_\al,\nu$ given in (\ref{Nu}). Then $\A:D(\A)\subset L^2_{\beta}(0,1)\rightarrow L^2_{\beta}(0,1)$ is a self-adjoint operator. Furthermore, the family $\{\Fi_k\}_ {k\geq 0}$ given by
		\begin{equation}\label{Phik}
			\begin{array}{rcl}
				\Fi_0(x) & := & \sqrt{2-\al+2\sqrt{\mu(\al+\beta)-\mu}}\,\,x^{(1-\al-\beta)/2+\sqrt{\mu(\al+\beta)-\mu}},\\
				\Fi_k(x)& := &\dfrac{\sqrt{2\kappa_\al}}{|J_\nu(j_{\nu+1,k})|}x^{(1-\al-\beta)/2}J_\nu\pts{j_{\nu+1,k}x^{\kappa_\al}},\quad k\geq 1,
			\end{array}
		\end{equation}
		is an orthonormal basis for $L^2_{\beta}(0,1)$ such that
		\begin{equation*}\label{lambdak}
			\mathcal{A} \Fi_k=\lambda_k \Fi_k, \quad k\geq 0,
		\end{equation*}
		where $\lambda_{0}:=0$ and $\lambda_{k}:=\kappa^{2}_{\al}(j_{\nu+1,k})^2,\, k\geq 1$.
	\end{proposition}
	\begin{proof} 
		Since $1/p,q,w\in L^1(1/2,1)$ we have that $x=1$ is a regular point.\\
		
		\textit{Case i)} Assume $\sqrt{\mu(\al+\beta)-\mu}< \kappa_\al$. \\
		
		First, we will build a (BC) basis $\{y_{0},z_{0}\}$ at $x=0$ and a (BC) basis $\{y_{1}, z_{1}\}$ at $x=1$, see \cite[Definition 10.4.3]{Zettl}.\\
		
		Consider the functions given by
		\begin{equation}\label{BCleft}
			y_{0}(x):=x^{(1-\al-\beta)/2+\sqrt{\mu(\al+\beta)-\mu}},\,\, z_{0}(x):=\dfrac{x^{(1-\al-\beta)/2-\sqrt{\mu(\al+\beta)-\mu}}}{2\sqrt{\mu(\al+\beta)-\mu}},\quad x\in (0,1).
		\end{equation}
		Notice the assumption implies that $y_{0}, z_{0}\in D_{\max}$. Clearly, $[z_{0}, y_{0}](0)=1$, thus $\{y_{0}, z_{0}\}$ is a (BC) basis at $x=0$.\\ 
		
		Since $y_{0}, z_{0}\in L^2_\beta(0,1)$ are linearly independent solutions of $Mu=0u$ it follows that $x=0$ is limit-circle (LC), see \cite[Definition 7.3.1, Theorem 7.2.2]{Zettl}.\\
		
		Consider also the functions given by
		\begin{equation*}\label{BCright}
			y_{1}(x):=-x^{(1-\al-\beta)/2+\sqrt{\mu(\al+\beta)-\mu}},\,\,z_{1}(x):=\dfrac{x^{(1-\al-\beta)/2+\sqrt{\mu(\al+\beta)-\mu}}-x^{(1-\al-\beta)/2-\sqrt{\mu(\al+\beta)-\mu}}}{2\sqrt{\mu(\al+\beta)-\mu}},
		\end{equation*} 
		Since $y_{1}, z_{1}\in D_{\max}$ and $[z_{1}, y_{1}](1)=1$, it follows that $\{y_{1}, z_{1}\}$ is a (BC) basis at $x=1$.\\
		
		Now, we fix $c,d\in (0,1)$ with $c<d$. From the Patching Lemma, Lemma 10.4.1 in \cite{Zettl}, there exist functions $g_1, g_2\in D_{\max}$ such that
		\begin{equation*}
			\begin{cases}g_1(c)=y_0(c), & g_1(d)=y_1(d), \\ (p g_1^{\prime})(c)=(p y_0^{\prime})(c), & (p g_1^{\prime})(d)=(p y_1^{\prime})(d),
			\end{cases}
			\,\,\,\,\,
			\begin{cases}g_2(c)=z_0(c), & g_2(d)=z_1(d), \\ (p g_2^{\prime})(c)=(p z_0^{\prime})(c), & (p g_2^{\prime})(d)=(p z_1^{\prime})(d).
			\end{cases}
		\end{equation*}
		Thus, the pair $\{y_+,y_-\}$ is a (BC) basis on $(0,1)$, see \cite[Definition 10.4.3]{Zettl}, where
		\begin{equation*}\label{BC}
			y_{+}(x):=\left\{\begin{array}{lll}
				y_0(x) & \text { if } & x \in(0, c), \\
				g_1(x) & \text { if } & x \in[c, d], \\
				y_1(x) & \text { if } & x \in(d, 1),
			\end{array}\right.\,\,\,\,\,\,
			y_{-}(x):=\left\{\begin{array}{lll}
				z_0(x) & \text { if } & x \in(0, c), \\
				g_2(x) & \text { if } & x \in[c, d], \\
				z_1(x) & \text { if } & x \in(d, 1).
			\end{array}\right.
		\end{equation*}
		
		The matrices
		$$A=\begin{pmatrix}1 & 0\\ 0 & 0\end{pmatrix}\quad\text{and}\quad B=\begin{pmatrix}0 & 0\\ 1 & 0\end{pmatrix}$$
		satisfy the hypothesis in \cite[Theorem 10.4.2]{Zettl}, then
		$$
		\begin{aligned}
			D(\A): & =\left\{u \in D_{\max }: A \left(\begin{array}{l}
				\left[u, y_{+}\right](0) \\
				{\left[u, y_{-}\right](0)}
			\end{array}\right)+B\left(\begin{array}{l}
				{\left[u, y_{+}\right](1)} \\
				{\left[u, y_{-}\right](1)}
			\end{array}\right)=\left(\begin{array}{l}
				0 \\
				0
			\end{array}\right)\right\} \\
			& =\left\{u \in D_{\max }:\left[u, y_{+}\right](0)=\left[u, y_{+}\right](1)=0\right\}=\left\{u \in D_{\max }:\left[u, y_{+}\right](0)=(au+u_x)(1)=0\right\}
			%& =\left\{u \in D_{\max }:\left[u, y_{+}\right](0)=0,\,\,\,\, (au+u')(1)=0\right\}
		\end{aligned}
		$$ 
		is a self-adjoint domain, therefore the operator $\mathcal{A}:D(\A)\subset L^2_\beta(0,1)\rightarrow L^2_\beta(0,1)$ is self-adjoint.\\
		
		Finally, we have that 
		\[\left[u, y_{+}\right](0)=\lim_{x\rightarrow 0^+}[u, y_{0}](x)=\lim_{x\rightarrow 0^+}\lvs{\frac{u}{z_0}(x)[z_0,y_0](x)+[u,z_0](x)\frac{y_0}{z_0}(x)}=\lim_{x\rightarrow 0^+}\frac{u}{z_0}(x),\]
		because  $[z_0,y_0](0)=1$, $[u,z_0](0)$ is finite (see \cite[Lemma 10.2.3]{Zettl}), and $\lim_{x\rightarrow 0^+}y_0/z_0(x)=0$. Hence, the result follows.\\
		
		\textit{Case ii)} Assume $\sqrt{\mu(\al+\beta)-\mu}\geq\kappa_\al$.\\
		
		The assumption implies that $z_{0}\notin L^2_{\beta}(0,1)$, then $x=0$ is limit point (LP). Theorem 10.4.4 in \cite{Zettl} with $A_1=a, A_2=1$ implies that $D(\A)=\{u\in D_{\max}\,| (au+u_{x})(1)=0\}$ is a self-adjoint domain.\\
		
		This concludes the first part of the proof.\\
		
		Clearly, $\Phi_k\in C^\infty(0,1)$ and (\ref{asincero}) implies that $\Phi_k\in L^2_\beta(0,1)$ for all $k\geq 0$. Moreover,
		\[\lim_{x\rightarrow 0^+}x^{(\al+\beta-1)/2+\sqrt{\mu(\al+\beta)-\mu}}\Phi_k(x)=C_{\al,\beta,\mu}\lim_{x\rightarrow 0^+}x^{2\sqrt{\mu(\al+\beta)-\mu}}=0,\quad k\geq 0.\]
		
		By using (\ref{recur}) we obtain
		\begin{eqnarray*}
			\dfrac{|J_\nu(j_{\nu+1,k})|}{\sqrt{2\kappa_\al}}\Fi'_k(1)&=& \frac{1-\al-\beta}{2}J_\nu(j_{\nu+1,k})+\kappa_\al j_{\nu+1,k}J'_\nu(j_{\nu+1,k})\\
			&=&\pts{\frac{1-\al-\beta}{2}+\kappa_\al \nu}J_\nu(j_{\nu+1,k})=-a\dfrac{|J_\nu(j_{\nu+1,k})|}{\sqrt{2\kappa_\al}}\Fi_k(1),
		\end{eqnarray*}
		therefore $(a\Fi_k+\Fi'_k)(1)=0$ for all $k\geq 1$. Clearly $(a\Fi_0+\Fi'_0)(1)=0$. Therefore, $\Phi_k\in D(\mathcal{A})$ for all $k\geq 0$.\\
		
		We set $v(x)=x^bJ_\nu(cx^r)$ with $r,c>0$ and $b\in\R$. In the proof of Proposition 11 in \cite{GaloLopez} was shown that
		\[x^{2-2r}\frac{d^2v}{dx^2}+(1-2b)x^{1-2r}\frac{dv}{dx}+(b^2-r^2\nu^2 )x^{-2r}v=-r^2c^2v.\]
		By taking $r=\kappa_\al, b=(1-\al-\beta)/2$, and $c=j_{\nu+1,k}$, we get $\A\Phi_k=\lambda_k \Phi_k$ for all $k\geq 1$. Clearly, $\A\Phi_0=0$.\\
		The result follows by Proposition \ref{basis}.
	\end{proof}
	
	\begin{remark}\label{boundcon}
		If $\sqrt{\mu(\al+\beta)-\mu}\geq\kappa_\al$, from Lemma 10.4.1(b) in \cite{Zettl} we have that $[u,y_0](0)=0$ for all $u \in D(\mathcal{A})$. When $\sqrt{\mu(\al+\beta)-\mu}<\kappa_\al$, in the proof of the last proposition was shown that $[u,y_0](0)=0$ for all $u \in D(\mathcal{A})$, where $y_0$ is given in (\ref{BCleft}).
	\end{remark}
	
	\begin{remark}
		The family $\{\Theta_k\}_{k\geq 0}$ given in (\ref{Dini}) is the so-called Fourier-Dini basis for $L^2(0,1)$.
	\end{remark}
	
	Then $(\A,D(\A))$ is the infinitesimal generator of a diagonalizable self-adjoint semigroup in $L^2_{\beta}(0,1)$. Thus, we can consider interpolation spaces for the initial data. For any $s\geq 0$, we define
	\[\mathcal{H}^{s}=\mathcal{H}^{s}(0,1):=D(\mathcal{A}^{s/2})=\left\{u=\sum_{k=0}^\infty a_{k} \Phi_{k}:\|u\|_{\mathcal{H}^{s}}^{2}=|a_{0}|^2+\sum_{k=1}^\infty |a_{k}|^{2} \lambda_{k}^{s}<\infty\right\},\]
	and we also consider the corresponding dual spaces
	\[\mathcal{H}^{-s}:=\left[\mathcal{H}^{s}(0,1)\right]^{\prime}.\]
	It is well known that $\mathcal{H}^{-s}$ is the dual space of $\mathcal{H}^{s}$ with respect to the pivot space $L^2_\beta(0,1)$, i.e
	\[\mathcal{H}^s\hookrightarrow \mathcal{H}^0=L^2_{\beta}(0,1)=\left(L^2_{\beta}(0,1)\right)'\hookrightarrow \mathcal{H}^{-s},\quad s>0. \]
	Equivalently, $\mathcal{H}^{-s}$ is the completion of $L^2_\beta(0,1)$ with respect to the norm
	\[\|u\|^2_{-s}:=|\langle u,\Phi_{0}\rangle_{\beta}|^{2}+\sum_{k=1}^{\infty}\lambda_k^{-s}|\langle u,\Phi_k\rangle_{\beta}|^2.\]
	It is well known that the linear mapping given by
	\[S(t)u_0=\sum_{k=0}^\infty \textrm{e}^{-\lambda_k t}a_k\Fi_k\quad\text{if}\quad u_0=\sum_{k=0}^\infty a_{k} \Phi_{k}\in \mathcal{H}^s,\]
	defines a self-adjoint semigroup $\{S(t)\}_{t\geq 0}$ in $\mathcal{H}^s $ for all $s\in\R$.\\
	
	For $\delta\in\R$ and a function $h:(0,1)\rightarrow \R$ we introduce the notion of $\delta$-generalized limit of $h$ at $x=0$ as follows
	\[\OO_\delta (h):=\lim_{x\rightarrow 0^+} x^{\delta}h(x).\]
	
	\textbf{Notation} Let $t>0$ fixed. If $z\in\HH^s$ then $S(t)z\in\HH^s$, so we write $\lim_{x \rightarrow 1^-} S(t)z$ instead of $\lim_{x \rightarrow 1^-} (S(t)z)(x)$.\\
	
	\subsection{Notion of weak solutions for both systems}
	Now we consider a convenient definition of a weak solution for the system (\ref{problem1Left}). Let $\tau>0$ be fixed. We multiply the equation in (\ref{problem1Left}) by $x^\beta\varphi(x,t)=x^\beta S(\tau-t)z^{\tau}$, $0\leq t\leq\tau$, integrate by parts (formally) and by using the boundary conditions for $u,\varphi$, see Remark \ref{boundcon}, we get
	\begin{eqnarray*}
		\langle u(\tau), z^{\tau}\rangle_\beta-\langle u_{0}, S(\tau)z^\tau\rangle_\beta&=&\int_0^T[u(\cdot,t),S(\tau-t)z^\tau](0)\mathrm{d}t\\
		&=&\int_0^T[u(\cdot,t),x^{-a}](0)\mathcal{O}_a(S(\tau-t)z^\tau)\mathrm{d}t
		=\int_0^Tf(t)\mathcal{O}_a(S(\tau-t)z^\tau)\mathrm{d}t.
	\end{eqnarray*}
	\begin{definition}
		Let $T>0$, $0\leq\al<2$, $\beta\in \R$, $\mu<\mu(\alpha+\beta)$, and $a$ given by (\ref{a_const}). Let $f \in L^2(0,T)$ and $u_0\in \HH^{-s}$ for some $s > 0$. A weak solution of (\ref{problem1Left}) is a function $u \in C^0([0,T];\HH^{-s})$ such that for every $\tau \in (0,T]$ and for every $z^\tau \in \HH^s$ we have
		\begin{equation}\label{weaksolLeft}
			\left\langle u(\tau), z^{\tau}\right\rangle_{\mathcal{H}^{-s}, \mathcal{H}^{s}} = \left\langle u_{0}, S(\tau)z^\tau\right\rangle_{\mathcal{H}^{-s}, \mathcal{H}^{s}} +\int_{0}^{\tau} f(t) \mathcal{O}_{a}\left(S(\tau-t) z^{\tau}\right) \mathrm{d} t.
		\end{equation}
	\end{definition}
	
	The next result shows the existence of weak solutions for the system (\ref{problem1Left}) under suitable conditions on the parameters $\alpha,\beta,\mu,$ and $s$, its proof is similar to the proof of Proposition 10 in \cite{GaloLopez}.
	\begin{proposition}\label{continuityL}
		Let $T>0$, $0\leq\al<2$, $\beta \in \R$, $\mu<\mu(\alpha+\beta)$, $a$ given in (\ref{a_const}). Let $f \in L^2(0,T)$ and $u_0\in \HH^{-s}$ such that $s>\nu$, with $\nu$ given in (\ref{Nu}). Then, formula (\ref{weaksolLeft}) defines for each $\tau \in [0, T ]$ a unique element $u(\tau) \in \HH^{-s}$ that can be written as
		\[
		u(\tau) = S(\tau) u_{0} + B(\tau) f, \quad \tau \in(0, T],\]
		where $B(\tau)$ is the strongly continuous family of bounded operators $B(\tau): L^{2}(0,T) \rightarrow \mathcal{H}^{-s}$ given by
		\[\left\langle B(\tau) f, z^{\tau}\right\rangle_{\mathcal{H}^{-s}, \mathcal{H}^{s}}=\int_{0}^{\tau} f(t) \mathcal{O}_{a}\left(S(\tau-t) z^{\tau}\right) \mathrm{d} t, \quad \text{for all  }z^{\tau} \in \mathcal{H}^{s} .\]
		Furthermore, the unique weak solution $u$ on $[0, T]$ to (\ref{problem1Left}) (in the sense of (\ref{weaksolLeft})) belongs to ${C}^{0}\left([0, T] ; \mathcal{H}^{-s}\right)$ and fulfills
		\[
		\|u\|_{L^{\infty}\left([0, T] ; \mathcal{H}^{-s}\right)} \leq C\left(\left\|u_{0}\right\|_{\mathcal{H}^{-s}}+\|f\|_{L^{2}(0, T)}\right).
		\]
	\end{proposition}
	\begin{proof}
		Fix $\tau>0$. Let $u(\tau)\in H^{-s}$ be determined by the condition (\ref{weaksolLeft}), hence
		$$
		u(\tau) - S(\tau) u_{0}=\zeta(\tau)f,
		$$
		where
		$$
		\left\langle\zeta(\tau)f, z^{\tau}\right\rangle_{\mathcal{H}^{-s}, \mathcal{H}^{s}}=\int_{0}^{\tau} f(t) \mathcal{O}_{a}\left(S(\tau-t) z^{\tau}\right) \mathrm{d} t, \quad \text{for all  } z^{\tau} \in \mathcal{H}^{s}.
		$$
		We claim that $\zeta(\tau)$ is a bounded operator from $L^{2}(0, T)$ into $\mathcal{H}^{-s}$: consider $z^{\tau} \in \mathcal{H}^{s}$ given by 
		\begin{equation}\label{finaldata}
			z^{\tau}=\sum_{k=0}^\infty b_{k} \Phi_{k},
		\end{equation}
		therefore
		$$
		S(\tau-t) z^{\tau}=\sum_{k=0}^\infty \mathrm{e}^{\lambda_{k}(t-\tau)} b_{k} \Phi_{k}, \quad \text{for all } t \in[0, \tau].
		$$
		By using Lemma \ref{reduce} and (\ref{lim2}) we obtain that there exists a constant $C=C(\al,\beta,\mu)>0$ such that
		\[|\mathcal{O}_{a}\left(\Phi_k\right)|\leq C |j_{\nu+1,k}|^{\nu+1/2},\quad k\geq 1,\]
		hence (\ref{below}) implies that there exists a constant $C=C(\al,\beta,\mu,\tau)>0$ such that
		\begin{eqnarray*}
			\pts{\int_{0}^{\tau}\left|\mathcal{O}_{a}\left(S(\tau-t) z^{\tau}\right)\right|^{2} \mathrm{~d} t}^{1/2}&\leq & \sum_{k=0}^\infty 
			|b_k| |\OO_{a}(\Fi_k)| \pts{\int_0^\tau \mathrm{e}^{2\lambda_{k}(t-\tau)} \mathrm{~d} t}^{1/2} \\
			&\leq& C\pts{\tau^{1/2}|b_{0}|+\pts{\sum_{k=1}^{\infty}|b_{k}|^2\lambda_{k}^{s}}^{1/2}\pts{\sum_{k=1}^\infty |\lambda_k|^{\nu-1/2-s}\pts{1-\mathrm{e}^{-2\lambda_k \tau}}}^{1/2}}\\
			&\leq & C\pts{\tau^{1/2}|b_{0}|+\pts{\sum_{k=1}^{\infty}|b_{k}|^2\lambda_{k}^{s}}^{1/2}\pts{\sum_{k=1}^\infty\frac{1}{k^{2(s-\nu+1/2)}}}^{1/2}}\\
			&\leq& C\left\|z^{\tau}\right\|_{\mathcal{H}^{s}}.
		\end{eqnarray*}
		Therefore $\|\zeta(\tau) f\|_{\mathcal{H}^{-s}}\leq C\|f\|_{L^2(0,T)}$ for all $f\in L^2(0,T)$, $\tau\in (0,T]$.\\
		
		Finally, we fix $f\in L^2(0,T)$ and show that the mapping $\tau\mapsto \zeta(\tau) f$ is right-continuous on $[0,T)$. Let $h>0$ small enough and $z\in \mathcal{H}^s$ given as in (\ref{finaldata}). Thus, proceeding as in the last inequalities, we have
		\begin{equation*}
			\begin{array}{l}
				\d|\left\langle\zeta(\tau+h)f-\zeta(\tau)f, z\right\rangle_{\mathcal{H}^{-s}, \mathcal{H}^{s}}|\\
				\leq
				\d C\|f\|_{L^2(0,T)}\pts{|b_{0}|h+\pts{\sum_{k=1}^{\infty}|b_{k}|^2\lambda_{k}^{s}}^{1/2}\cts{\pts{\sum_{k=1}^\infty\frac{I(\tau,k,h)}{k^{2(s-\nu+1/2)}}}^{1/2}+\pts{\sum_{k=1}^\infty\frac{1-\mathrm{e}^{-2\lambda_k h}}{k^{2(s-\nu+1/2)}}}^{1/2}}},
			\end{array}
		\end{equation*}
		where 
		\begin{equation}\label{Ianu}
			I(\tau,k,h)=\lambda_k\int_0^\tau\pts{\mathrm{e}^{\lambda_k(t-\tau-h)}-\mathrm{e}^{\lambda_k(t-\tau)}}^2\mathrm{~d} t
			= \frac{1}{2}(1-\mathrm{e}^{-\lambda_k h})^2(1-\mathrm{e}^{-2\lambda_k \tau})\rightarrow 0\quad\text{as}\quad h\rightarrow 0^+.
		\end{equation}
		Since $0\leq I(\tau,k,h)\leq 1/2$ uniformly for $\tau, h>0$, $k\geq 1$, the result follows by the dominated convergence theorem.
	\end{proof}
	
	\begin{remark}
		In the following section, we will consider initial conditions in $L^2_\beta(0,1)$. Notice that $L^2_\beta(0,1)\subset H^{-\nu-\delta}$ for all $\delta>0$, and we can apply Proposition \ref{continuityL} with $s=\nu+\delta$, $\delta>0$, then the corresponding solutions will be in $C^0([0, T ], H^{-\nu-\delta})$.\\
	\end{remark}
	
	As before, we introduce a suitable definition of a weak solution for the system (\ref{problem1Right}).
	\begin{definition}
		Let $T>0$, $\beta\in \R$, $0\leq\al<2$, $\mu<\mu(\alpha+\beta)$ and $a$ given in (\ref{a_const}). Let $f \in L^2(0,T)$ and $u_0\in L_{\beta}^2(0,1)$. A weak solution of (\ref{problem1Right}) is a function $u \in C^0\pts{[0,T];L_{\beta}^2(0,1)}$ such that for every $\tau \in (0,T]$ and for
		every $z^\tau \in L^2_{\beta}(0,1)$ we have
		\begin{equation}\label{weaksolP1R}
			\left\langle u(\tau), z^{\tau}\right\rangle_{\beta} = \left\langle u_{0}, S(\tau)z^\tau\right\rangle_{\beta} + \int_{0}^{\tau} f(t) \lim_{x\to 1^{-}}S(\tau-t) z^{\tau} \mathrm{d} t.
		\end{equation}
	\end{definition}
	The next result shows the existence of weak solutions for the system (\ref{problem1Right}) under certain conditions on the parameters $\alpha,\beta,\mu$ and $a$, its proof is similar to the proof of Proposition 18 in \cite{GaloLopez2}.
	\begin{proposition}\label{continuity1R}
		Let $T>0$, $\beta\in \R$, $0\leq\al<2$, $\mu<\mu(\alpha+\beta)$ and $a$ given in (\ref{a_const}). Let $f \in L^2(0,T)$ and $u_0\in L_{\beta}^2(0,1)$. Then, formula (\ref{weaksolP1R}) defines for each $\tau \in [0, T ]$ a unique element $u(\tau) \in L_{\beta}^2(0,1)$ that can be written as
		\[
		u(\tau) - S(\tau) u_{0} = \mathcal{B}(\tau) f, \quad \tau \in(0, T],\]
		where $B(\tau)$ is the strongly continuous family of bounded operators $\mathcal{B}(\tau): L^{2}(0,T) \rightarrow L_{\beta}^2(0,1)$ given by
		\[\left\langle \mathcal{B}(\tau) f, z^{\tau}\right\rangle_{\beta}=\int_{0}^{\tau} f(t) \lim_{x\to 1^{-}}S(\tau-t) z^{\tau} \mathrm{d} t, \quad \text{for all  }z^{\tau} \in L_{\beta}^2(0,1) .\]
		Furthermore, the unique weak solution $u$ on $[0, T]$ to (\ref{problem1Right}) (in the sense of (\ref{weaksolP1R})) belongs to ${C}^{0}\left([0, T] ; L_{\beta}^2(0,1)\right)$ and fulfills
		\[
		\|u\|_{L^{\infty}\left([0, T] ; L_{\beta}^2(0,1)\right)} \leq C\left(\left\|u_{0}\right\|_{\beta}+\|f\|_{L^{2}(0, T)}\right).
		\]
	\end{proposition}
	\begin{proof}
		Fix $\tau>0$. Let $u(\tau)\in L_{\beta}^2(0,1)$ be determined by the condition (\ref{weaksolP1R}), hence
		$$
		u(\tau) - S(\tau) u_{0}=\zeta(\tau)f,
		$$
		where
		$$
		\left\langle\zeta(\tau)f, z^{\tau}\right\rangle_{\beta}=\int_{0}^{\tau} f(t) \lim_{x\rightarrow 1^- }S(\tau-t) z^{\tau}\mathrm{d} t \quad \text{for all  } z^{\tau} \in L_{\beta}^2(0,1).
		$$
		Let $z^{\tau} \in L_{\beta}^2(0,1)$ written as
		\begin{equation}\label{finaldata1R}
			z^{\tau}=\sum_{k=0}^\infty b_{k} \Phi_{k},
		\end{equation}
		therefore
		$$
		\lim_{x\rightarrow 1^- }S(\tau-t) z^{\tau}=\sum_{k=0}^\infty \mathrm{e}^{\lambda_{k}(t-\tau)} b_{k} \Phi_{k}(1) \quad \text{for all } t \in[0, \tau].
		$$
		By (\ref{Phik}) we get
		\begin{equation}\label{Phik_1}
			\left|\Phi_{0}(1)\right|=\sqrt{2-\al+2\sqrt{\mu(\al+\beta)-\mu}},\,\,\quad\quad\left|\Phi_{k}(1)\right| = \sqrt{2\kappa_{\al}}, \quad k\geq 1,
		\end{equation}
		hence there exists a constant $C=C(\al,\beta,\mu,\tau)>0$ such that
		\begin{eqnarray*}
			\pts{\int_{0}^{\tau}\left|\lim_{x\rightarrow 1^- }S(\tau-t) z^{\tau}\right|^{2} \mathrm{~d} t}^{1/2}&\leq & \sum_{k=0}^\infty 
			|b_k| |\Fi_k(1)| \pts{\int_0^\tau \mathrm{e}^{2\lambda_{k}(t-\tau)} \mathrm{~d} t}^{1/2} \\
			&\leq & C\|z^{\tau}\|_{\beta}\pts{\sum_{k=0}^\infty \int_0^\tau \mathrm{e}^{2\lambda_{k}(t-\tau)} \mathrm{~d} t}^{1/2}
			=C\|z^{\tau}\|_{\beta}\pts{\tau+\sum_{k=1}^\infty \frac{1-\mathrm{e}^{-2\lambda_k \tau}}{2\lambda_k}}^{1/2}\\
			&\leq &C\|z^{\tau}\|_{\beta}\pts{\tau+\sum_{k=1}^\infty \frac{1}{k^2}}^{1/2}.
		\end{eqnarray*}
		Therefore $\|\zeta(\tau) f\|_{\beta}\leq C\|f\|_{L^2(0,T)}$ for all $f\in L^2(0,T)$, $\tau\in (0,T]$.\\
		
		Finally, we fix $f\in L^2(0,T)$ and show that the mapping $\tau\mapsto \zeta(\tau) f$ is right-continuous on $[0,T)$. Let $h>0$ small enough and $z\in L^2_\beta(0,1)$ given as in (\ref{finaldata1R}). Then we have
		\begin{eqnarray*}
			|\left\langle\zeta(\tau+h)f-\zeta(\tau)f, z\right\rangle_{\beta}|&\leq&\int_{0}^{\tau} |f(t)|\left|\lim_{x\rightarrow 1^- }(S(\tau+h-t)-S(\tau-t)) z\right| \mathrm{d} t \\
			&&+ \int_{\tau}^{\tau+h} |f(t)|\left|\lim_{x\rightarrow 1^- }S(\tau+h-t) z\right| \mathrm{d} t\\
			&\leq& \d C\|z^{\tau}\|_{\beta}\|f\|_{L^2(0,T)}\cts{\pts{\sum_{k=1}^\infty\frac{I(\tau,k,h)}{k^{2}}}^{1/2}+\pts{h+\sum_{k=1}^\infty\frac{1-\mathrm{e}^{-2\lambda_k h}}{k^{2}}}^{1/2}},
		\end{eqnarray*}
		where $I(\tau,k,h)\rightarrow 0$ as $h\rightarrow 0^+$, see (\ref{Ianu}).
	\end{proof}
	%%%%%%%%%%%%%%%%%%%%%%%%%%%%%%%%%%%%%%%%%%
	\section{Control at the left endpoint}\label{leftend}
	
	\subsection{Upper estimate of the cost of the null controllability}\label{P1Control}
	Here we use the moment method, introduced by Fattorini \& Russell in \cite{Fatorini}, to prove the null controllability of the system (\ref{problem1Left}). The first step is to construct a biorthogonal family $\d \lvs{\psi_k}_{k\geq 0}\subset L^2(0,T)$ to the family of exponential functions $\lvs{\mathrm{e}^{-\lambda_{k}(T-t)}}_{k\geq 0}$ on $[0, T]$, i.e that satisfies
	$$
	\int_{0}^{T}\psi_k(t) \mathrm{e}^{-\lambda_{l}(T-t)} dt = \delta_{kl},\quad\text{for all}\quad k,l\geq 0.
	$$
	This construction will help us to get an upper bound for the cost of the null controllability of the system (\ref{problem1Left}).\\
	
	Assume that for each  $k\geq 0$ there exists an entire function $F_k$ of exponential type $T/2$ such that $F_k(x)\in L^2(\R)$, and
	\begin{equation}\label{krone}
		F_{k}(i\lambda_{l})=\delta_{kl},\quad \text{for all}\quad k,l\geq 0.
	\end{equation}
	The $L^2$-version of the Paley-Wiener theorem implies that there exists $\eta_k\in L^2(\R)$ with support in $[-T/2,T/2]$ such that $F_k(z)$ is the analytic extension of the Fourier transform of $\eta_k$. Then we have that 
	\begin{equation}\label{psieta}
		\psi_k(t):=\mathrm{e}^{\lambda_k T/2}\eta_k(t-T/2),\quad t\in[0,T],\,k\geq0,
	\end{equation}
	is the family we are looking for.\\
	%\[\delta_{kl}=\mathrm{e}^{(\lambda_k-\lambda_l)T/2}F_{k}(i\lambda_{l})=\mathrm{e}^{(\lambda_k-\lambda_l)T/2}\int_{-\frac{T}{2}}^{\frac{T}{2}} \eta_{k}(t) \mathrm{e}^{\lambda_{l}t} \mathrm{d}t=\int_{0}^{T}\psi_k(t) \mathrm{e}^{-\lambda_{l}(T-t)} \mathrm{d}t\quad\text{for all}\quad k,l\geq 1.\]
	
	Now, we proceed to construct the family $F_k$, $k\geq 0$. Consider the Weierstrass infinite product
	\begin{equation}\label{weier}
		\Lambda(z):=z\prod_{k=1}^{\infty}\pts{1+\dfrac{iz}{(\kappa_\al j_{\nu+1, k})^2}}.
	\end{equation}
	From (\ref{asint}) we have that $j_{\nu+1, k}=O(k)$ for $k$ large, thus the infinite product converges absolutely in $\C$. Hence $\Lambda(z)$ is an entire function with simple zeros at $i\lambda_k$, $k\geq 0$.\\
	
	From \cite[Chap. XV, p. 498, eq. (3)]{Watson}, we have for $\nu>-1$ that
	\begin{equation}\label{explicit}
		\Lambda(z)=z\Gamma(\nu+2)\pts{\dfrac{2\kappa_\al}{\sqrt{-iz}}}^{\nu+1}J_{\nu+1}\pts{\dfrac{\sqrt{-iz}}{\kappa_\al}}.
	\end{equation}
	
	In \cite{GaloLopez} was proved that
	\[|J_\nu(z)|\leq \frac{|z|^\nu \mathrm{e}^{|\Im(z)|}}{2^\nu\Gamma\left(\nu+1\right)},\quad z\in\C.\]
	Therefore,
	\[|\Lambda(z)|\leq |z|\exp\pts{\frac{|\Im(\sqrt{-iz})|}{\kappa_\alpha}},\quad z\in\C.\]
	In particular,
	\begin{equation}\label{besst}
		|\Lambda(z)|\leq|z|\exp\pts{\frac{|z|^{1/2}}{\kappa_\al}},\quad z\in\C,\quad |\Lambda(x)|\leq|x|\exp\pts{\frac{|x|^{1/2}}{\sqrt{2}\kappa_\al}},\quad x\in\R.
	\end{equation}
	It follows that 
	\begin{equation}\label{PsiFunction}
		\Psi_{k}(z):=\dfrac{\Lambda(z)}{\Lambda'(i\lambda_{k})(z-i\lambda_{k})},\quad k\geq 0,
	\end{equation}
	is a family of entire functions that satisfy (\ref{krone}). Since $\Psi_{k}(x)$ is not in $L^2(\R)$, we need to fix this by using a suitable ``complex multiplier", thus we follow the approach introduced in \cite{Tucsnak}.\\
	
	For $\theta,\omega>0$, we define
	$$
	\sigma_{\theta}(t):=\exp\pts{-\frac{\theta}{1-t^2}},\quad t\in(-1,1),
	$$
	and extended by $0$ outside of $(-1, 1)$. Clearly $\sigma_{\theta}$ is analytic on $(-1,1)$. Set $C_{\theta}^{-1}:=\int_{-1}^{1}\sigma_{\theta}(t)\mathrm{d}t$ and define
	\begin{equation}\label{Hfunction}
		H_{\omega,\theta}(z)=C_{\theta}\int_{-1}^{1}\sigma_{\theta}(t)\exp\pts{-i\omega tz}\mathrm{d}t.
	\end{equation}
	$H_{\omega,\theta}(z)$ is an entire function, and the next result provides additional properties of $H_{\omega,\theta}(z)$.\\
	\begin{lemma}
		The function $H_{\omega,\theta}$ fulfills the following inequalities
		\begin{eqnarray}
			H_{\omega,\theta}(ix)&\geq &\frac{\exp\pts{\omega|x|/\pts{2\sqrt{\theta+1}}}}{11\sqrt{\theta+1}},\quad x\in\R,\label{Hcot1}\\
			|H_{\omega,\theta}(z)|     &\leq & \exp\pts{\omega|\Im(z)|},\quad z\in\C,\label{Hcot2}\\
			|H_{\omega,\theta}(x)|     &\leq & \chi_{|x|\leq 1}(x)+c\sqrt{\theta+1}\sqrt{\omega\theta\abs{x}}\exp\pts{3\theta/4-\sqrt{\omega\theta\abs{x}}}\chi_{|x|> 1}(x),\quad x\in\R,\label{Hcot3}
		\end{eqnarray}
		where $c>0$ does not depend on $\omega$ and $\theta$.
	\end{lemma}
	We refer to \cite[pp. 85--86]{Tucsnak} for the details.\\
	
	For $k\geq 0$ consider the entire function $F_{k}$ given as
	\begin{equation}\label{Ffunction}
		F_{k}(z):=\Psi_{k}(z)\dfrac{H_{\omega,\theta}(z)}{H_{\omega,\theta}(i\lambda_{k})},\quad z\in\C.
	\end{equation}
	
	For $\delta\in(0,1)$ we set
	\begin{equation}\label{aConst}
		\omega:=\frac{T(1-\delta)}{2}>0,\quad \text{and}\quad \theta:=\dfrac{(1+\delta)^2}{\kappa_\al^2 T\pts{1-\delta}}>0.
	\end{equation}
	
	\begin{lemma}
		The function $F_{k}(z)$, $k\geq 0$, has the following properties:\\
		i) $F_{k}$ is of exponential type $T/2$.\\
		ii) $F_{k}\in L^1(\R)\cap L^2(\R)$.\\
		iii) $F_k$ satisfies (\ref{krone}).\\
		iv) Furthermore, there exists a constant $c>0$, independent of $T,\alpha$ and $\delta$, such that
		\begin{equation}\label{F0bound}
			\left\|F_{0}\right\|_{L^{1}(\R)} \leq C(T, \alpha,\delta)\quad\text{and}
		\end{equation}
		\begin{equation}\label{Fbound}
			\left\|F_{k}\right\|_{L^{1}(\R)} \leq \frac{C(T, \alpha,\delta)}{\lambda_k\left|\Lambda^{\prime}\left(i \lambda_{k}\right)\right|}  \exp\pts{-\frac{\omega\lambda_k}{2\sqrt{\theta+1}}},\quad k\geq 1,
		\end{equation} 	
		where
		\begin{equation}\label{upper}
			C(T, \alpha,\delta)=c\sqrt{\theta+1}\cts{\exp\pts{{\frac{1}{\sqrt{2}\kappa_\alpha}}}+\sqrt{\theta+1}\frac{\kappa_\alpha^2}{\delta^5}\exp\pts{\frac{3 \theta}{4}}}. 
		\end{equation}
	\end{lemma}
	\begin{proof}
		By using (\ref{besst}), (\ref{Hcot2}), (\ref{Ffunction}) and (\ref{aConst}) we get that $F_{k}$ is of exponential type $T/2$ for all $k\geq 0$. Moreover, by using (\ref{PsiFunction}) and (\ref{Ffunction}), we can see that $F_{k}$ fulfills (\ref{krone}).\\
		
		Now we use (\ref{besst}), (\ref{Hcot1}), (\ref{Hcot3}), (\ref{Ffunction}), and (\ref{aConst}) to get
		\begin{eqnarray*}
			\left|F_{k}(x)\right| &\leq& c\exp\pts{-\frac{\omega\lambda_k}{2\sqrt{\theta+1}}}\frac{\sqrt{\theta+1}|x|}{|\Lambda^{\prime}\left(i \lambda_{k}\right)||x^2+ \lambda_{k}^2|^{1/2}} |H_{\omega,\theta}(x)|\exp\pts{ \frac{|x|^{1/2}}{\sqrt{2}\kappa_\alpha}}
			\\
			&\leq& c \exp\pts{-\frac{\omega\lambda_k}{2\sqrt{\theta+1}}}\frac{\sqrt{\theta+1}}{\lambda_k|\Lambda^{\prime}\left(i \lambda_{k}\right)|}\cts{\mathrm{e}^{\frac{1}{\sqrt{2}\kappa_\alpha}}\chi_{|x|\leq 1}(x)+\sqrt{\theta+1}\sqrt{\omega\theta}\abs{x}^{3/2}\exp\pts{\frac{3\theta}{4}- \frac{\delta|x|^{1/2}}{\sqrt{2}\kappa_\alpha}}\chi_{|x|> 1}(x)},
		\end{eqnarray*}	
		for all $k\geq 1$. Since the function on the right-hand side is rapidly decreasing in $\R$, we have $F_{k}\in L^1(\R)\cap L^2(\R)$. Finally, the change of variable $y=(\kappa_\alpha)^{-1}\delta|x|^{1/2}/\sqrt{2}$ implies (\ref{Fbound}).\\
		
		When $k=0$ we have
		\[|F_0(x)|\leq \exp\pts{\frac{|x|^{1/2}}{\sqrt{2}\kappa_\al}} |H_{\omega,\theta}(x)|\leq \mathrm{e}^{\frac{1}{\sqrt{2}\kappa_\alpha}}\chi_{|x|\leq 1}(x)+\sqrt{\theta+1}\sqrt{\omega\theta\abs{x}}\exp\pts{\frac{3\theta}{4}- \frac{\delta|x|^{1/2}}{\sqrt{2}\kappa_\alpha}}\chi_{|x|> 1}(x),\]
		then we integrate on $\R$ and the result follows.
	\end{proof}
	
	Since $\eta_k, F_k\in L^1(\R)$, the inverse Fourier theorem yields 
	\[\eta_k(t)=\frac{1}{2\pi}\int_{\R}\mathrm{e}^{it\tau}F_k(\tau)\mathrm{d}\tau,\quad t\in\R, k\geq 0,\]
	hence (\ref{psieta}) implies that $\psi_k\in C([0,T])$. From (\ref{F0bound}) and (\ref{Fbound}) we have $\|\psi_0\|_{\infty}\leq C(T, \alpha,\delta)$ and
	\begin{equation}\label{psiL1}
		\|\psi_k\|_{\infty}\leq \frac{C(T, \alpha,\delta)}{\lambda_k\left|\Lambda^{\prime}\left(i \lambda_{k}\right)\right|}  \exp\pts{\frac{T\lambda_k}{2}-\frac{\omega\lambda_k}{2\sqrt{\theta+1}}},\quad k\geq 1.	
	\end{equation}
	
	Now, we are ready to prove the null controllability of the system (\ref{problem1Left}). Let $u_{0}\in L^{2}_\beta(0,1)$. Then consider its (generalized) Fourier-Dini series with respect to the orthonormal basis $\{\Phi_{k}\}_{k\geq 0}$,
	\begin{equation}\label{uoSerie}
		u_{0}(x)=\sum_{k=0}^{\infty} b_{k} \Phi_{k}(x).
	\end{equation}
	We set
	\begin{equation}\label{fserie}
		f(t):=-\sum_{k=0}^{\infty}\frac{b_{k} \mathrm{e}^{-\lambda_{k} T}}{\mathcal{O}_{a}\left(\Phi_{k}\right)} \psi_{k}(t).
	\end{equation}
	Since $\{\psi_k\}_{k\geq 0}$ is biorthogonal to $\{\mathrm{e}^{-\lambda_k(T-t)}\}_{k\geq 0}$ we have
	$$
	\int_{0}^{T} f(t) \mathcal{O}_{a}\left(\Phi_{k}\right) \mathrm{e}^{-\lambda_{k}(T-t)} \mathrm{d} t = -b_{k} \mathrm{e}^{-\lambda_{k} T} = -\left\langle u_{0}, \mathrm{e}^{-\lambda_{k} T}\Phi_{k}\right\rangle_\beta = -\left\langle u_{0}, \mathrm{e}^{-\lambda_{k} T}\Phi_{k}\right\rangle_{\mathcal{H}^{-s}, \mathcal{H}^{s}},\quad k\geq 0.
	$$
	Let $u\in C([0,T];H^{-s})$ that satisfies (\ref{weaksolLeft}) for all $\tau\in (0,T]$, $z^\tau\in H^s$. In particular, for $\tau=T$ we take $z^T=\Phi_k$, $k\geq 0$, then the last equality implies that
	$$
	\left\langle u(\cdot, T), \Phi_{k}\right\rangle_{\mathcal{H}^{-s}, \mathcal{H}^{s}}=0\quad \text{for all}\quad k \geq 0,
	$$
	hence $u(\cdot, T)=0$.\\
	
	It just remains to estimate the norm of the control $f$. From (\ref{psiL1}) and (\ref{fserie})  we get
	\begin{equation}\label{finty}
		C(T, \alpha,\delta)^{-1}\|f\|_{\infty} \leq  \frac{|b_0|}{|\mathcal{O}_{a}(\Phi_{0})|}+\sum_{k=1}^{\infty} \frac{\left|b_{k}\right|}{\left|\mathcal{O}_{a}\left(\Phi_{k}\right)\right|} \frac{1}{\lambda_k\left|\Lambda^{\prime}\left(i \lambda_{k}\right)\right|} \exp\pts{-\frac{T\lambda_k}{2}-\frac{\omega\lambda_k}{2\sqrt{\theta+1}}}.%\exp \left(\frac{3 \theta}{4}-\frac{\lambda_{k} T}{2}\right).\frac{\kappa_\alpha^3}{\delta^3}
	\end{equation}
	
	From (\ref{weier}), (\ref{explicit}), and (\ref{recur2}) (with $\nu+1$ instead of $\nu$), we get that
	\begin{equation}\label{deriLambda}
		\Lambda^{\prime}\left(i \lambda_{k}\right)=i \lambda_{k}\frac{2^{\nu+1}\Gamma(\nu+2)}{(j_{\nu+1, k})^{\nu+2}} \frac{-i}{2\kappa_\al^{2}} J_{\nu+1}^{\prime}(j_{\nu+1, k})= \frac{2^{\nu}\Gamma(\nu+2)}{(j_{\nu+1, k})^{\nu}} J_{\nu}(j_{\nu+1, k}), \quad k\geq 1,
	\end{equation}
	and by using (\ref{lim2}) we get
	\begin{equation*}\label{estimate2}
		\d\left|\mathcal{O}_{a}\left(\Phi_{k}\right)\Lambda^{\prime}\left(i \lambda_{k}\right)\right|=\frac{\Gamma(\nu+2)}{\Gamma(\nu+1)} \sqrt{2\kappa_\al}=(\nu+1)\sqrt{2\kappa_\al},\quad k\geq 1. %\dfrac{[\mu(\alpha+\beta)+\kappa_\al\nu]j_{\nu, k}\sqrt{2\kappa_\al}}{2}.
	\end{equation*}
	From (\ref{finty}), (\ref{below}), and using that $\lambda_k\geq \lambda_1$, it follows that %Hence, combining (\ref{estimate1}) and (\ref{estimate2}) we deduce that
	\begin{eqnarray*}
		C(T, \alpha,\delta)^{-1}\|f\|_{\infty} &\leq& \frac{|b_0|}{|\mathcal{O}_{a}(\Phi_{0})|}+\frac{1}{\sqrt{2}(\nu+1)\kappa_\al^{5/2}} \exp\pts{-\frac{T\lambda_1}{2}-\frac{\omega\lambda_1}{2\sqrt{\theta+1}}}\sum_{k=1}^{\infty} \frac{|b_{k}| }{(j_{\nu+1, k})^2}\\
		&\leq& \frac{|b_0|}{|\mathcal{O}_{a}(\Phi_{0})|}+\frac{c}{(\nu+1)\kappa_\al^{5/2}} \exp\pts{-\frac{T\lambda_1}{2}-\frac{\omega\lambda_1}{2\sqrt{\theta+1}}}\pts{\sum_{k=1}^{\infty} |b_{k}|^2}^{1/2}.
	\end{eqnarray*}
	Using the expression of $\omega,\theta$ given in (\ref{aConst}) and the facts $\theta>0$, $\delta\in(0,1)$, and $0<\kappa_\alpha\leq 1$, we get that
	\[\theta\leq \frac{4}{(1-\delta)\kappa_\alpha^2 T },\quad \sqrt{\theta+1}\leq \frac{2(1+T)^{1/2}}{(1-\delta)^{1/2}\kappa_\alpha T^{1/2}},\quad \sqrt{\theta+1}\leq \theta+1,\]
	therefore
	\begin{equation}\label{reduc}
		\frac{\omega}{\sqrt{\theta+1}}\geq \frac{\kappa_\alpha(1-\delta)^{3/2}T^{3/2}}{4(1+T)^{1/2}}, \quad C(T, \alpha,\delta)\leq c\pts{1+\frac{1}{(1-\delta)\kappa_\alpha^2 T}}\cts{\exp\pts{\frac{1}{\sqrt{2}\kappa_\alpha}}+\frac{1}{\delta^5}\exp\pts{\frac{3}{(1-\delta)\kappa_\alpha^2 T}}}.
	\end{equation}
	By using the definition of $\lambda_{1}$, and setting $b_0=0$, we get the estimate for $\mathcal{K}_{\Phi_0^\perp}$.
	%%%%%%%%%%%%%%%%%%%%%%%%%%%%%%%%%%%%%
	%%%%%%%%%%%%%%%%%%%%%%%%%%%%%%
	\subsection{Lower estimate of the cost of the null controllability}
	In this section, we get a lower estimate of the cost $\mathcal{K}=\mathcal{K}_{L^2_\beta}(T,\alpha,\beta,\mu)$.\\
	We set
	\begin{equation}\label{first4} %x^{\pts{1-\al-\beta}/2}J_\nu\pts{j_{\nu,1}x^{\kappa_\al}}
		u_0(x):= \frac{|J_\nu(j_{\nu+1,1})|}{(2\kappa_{\al})^{1/2}}\Phi_{1}(x),\,x\in(0,1),
		\quad \text{hence}\quad\|u_0\|^2_\beta=\frac{|J_\nu(j_{\nu+1,1})|^2}{2\kappa_{\al}}.
	\end{equation}
	For $\varepsilon>0$ small enough, there exists $f\in U(\al,\beta,\mu,T,u_0)$ such that
	\begin{equation}\label{inicost4}
		u(\cdot,T)\equiv 0,\quad \text{and}\quad \|f\|_{L^2(0,T)}\leq (\mathcal{K}+\varepsilon)\|u_0\|_\beta.
	\end{equation}
	Then, in (\ref{weaksolLeft}) we set $\tau=T$ and  take $z^\tau=\Fi_k$, $k\geq 0$,  to obtain
	\begin{eqnarray*}
		\mathrm{e}^{-\lambda_k T}\left\langle u_{0},\Fi_k\right\rangle_\beta=\left\langle u_{0}, S(T)\Fi_k\right\rangle_{\mathcal{H}^{-s}, \mathcal{H}^{s}}&=&
		-\int_{0}^{T} f(t) \mathcal{O}_{a}\left(S(T-t) \Fi_k\right) \mathrm{d} t\\
		&=&
		-\mathrm{e}^{-\lambda_k T}\mathcal{O}_{a}\left(\Fi_k\right)\int_{0}^{T} f(t) \mathrm{e}^{\lambda_k t} \mathrm{d} t,
	\end{eqnarray*}
	from (\ref{first4}) and (\ref{lim2}) it follows that
	\begin{equation}\label{ortocont4}
		\int_{0}^{T} f(t) \mathrm{e}^{\lambda_k t} \mathrm{d} t= -\frac{2^{\nu}\Gamma(\nu+1)|J_\nu(j_{\nu+1,1})|^2}{2\kappa_{\al}(j_{\nu+1,1})^{\nu}}\delta_{1,k},\quad k\geq 0.
	\end{equation}
	
	Now consider the function $v: \mathbb{C} \rightarrow \mathbb{C}$ given by
	\begin{equation*}
		v(s):=\int_{-T / 2}^{T / 2} f\left(t+\frac{T}{2}\right) \mathrm{e}^{-ist} \mathrm{~d} t, \quad s \in \mathbb{C} .
	\end{equation*}
	Fubini and Morera's theorems imply that $v(s)$ is an entire function. Moreover, (\ref{ortocont4}) implies that
	\[v(i\lambda_k)=0\quad\text{for all }k\geq 0, k \neq 1, \quad \text{and}\quad v(i\lambda_1)=-\frac{2^{\nu}\Gamma(\nu+1)|J_\nu(j_{\nu+1,1})|^2}{2\kappa_{\al}(j_{\nu+1,1})^{\nu}}\mathrm{e}^{-\lambda_{1}T/2} .\]
	We also have that
	\begin{equation}\label{uve4}
		|v(s)| \leq \mathrm{e}^{T|\Im(s)|/2}\int_{0}^{T}|f(t)| \mathrm{d} t \leq (\mathcal{K}+\varepsilon) T^{1/2}\mathrm{e}^{T|\Im(s)|/2} \left\|u_{0}\right\|_{\beta}.
	\end{equation}
	Consider the entire function $F(z)$ given by
	\begin{equation}\label{entire4}
		F(s):=v\left(s-i \delta\right), \quad s \in \mathbb{C},
	\end{equation}
	for some $\delta>0$ that will be chosen later on. Clearly, 
	\[ F\left(b_{k}\right)=0, \quad k\geq 0, k \neq 1, \quad \text {where} \quad b_{k}:=i\left(\lambda_k+\delta\right),\quad k\geq 0,\quad\text{and}\]
	\begin{equation}\label{ena4}
		F(b_{1})=-\frac{2^{\nu}\Gamma(\nu+1)|J_\nu(j_{\nu+1,1})|^2}{2\kappa_{\al}(j_{\nu+1,1})^{\nu}}\mathrm{e}^{-\lambda_{1}T/2}.
	\end{equation}
	From (\ref{first4}), (\ref{uve4}) and (\ref{entire4}) we obtain
	\begin{equation}\label{logF4}
		\log |F(s)|\leq \frac{T}{2}|\Im(s)-\delta|+\log\pts{(\mathcal{K}+\varepsilon) T^{1 / 2}\frac{\abs{J_\nu\pts{j_{\nu+1,1}}}}{\pts{2\kappa_\al}^{1/2}}},\quad s\in\mathbb{C}.
	\end{equation}
	
	We recall the following representation theorem, see \cite[p. 56]{koosis}.
	\begin{theorem}\label{repreTheo} Let $g(z)$ be an entire function of exponential type and assume that
		$$
		\int_{-\infty}^{\infty} \frac{\log ^{+}|g(x)|}{1+x^{2}} \mathrm{d}x<\infty.
		$$
		Let $\left\{d_{\ell}\right\}_{\ell \geq 1}$ be the set of zeros of $g(z)$ in the upper half plane $\Im(z)>0$ (each zero being repeated as many times as its multiplicity). Then,
		$$
		\log |g(z)|=A \Im(z)+\sum_{\ell=1}^{\infty} \log \left|\frac{z-d_{\ell}}{z-\bar{d}_{\ell}}\right|+\frac{\Im(z)}{\pi} \int_{-\infty}^{\infty} \frac{\log |g(s)|}{|s-z|^{2}} \mathrm{d}s,\quad\Im(z)>0,
		$$
		where
		$$
		A=\limsup _{y \rightarrow\infty} \frac{\log |g(i y)|}{y} .
		$$
	\end{theorem}
	We apply the last result to the function $F(z)$ given in (\ref{entire4}). In this case, (\ref{uve4}) implies that $A\leq T/2$. Also notice that $\Im\left(b_{k}\right)>0$, $k\geq 0$, to get
	\begin{equation}\label{aprep4}
		\log \left|F\left(b_{1}\right)\right|\leq(\lambda_1+\delta)\frac{T}{2}+\sum_{k=0,k\neq 1}^{\infty} \log \left|\frac{b_{1}-b_{k}}{b_{1}-\bar{b}_{k}}\right|+\frac{\Im\left(b_{1}\right)}{\pi} \int_{-\infty}^{\infty} \frac{\log |F(s)|}{\left|s-b_{1}\right|^{2}} \mathrm{~d}s.
	\end{equation}
	By using the definition of the constants $b_k$'s we have
	\begin{eqnarray}
		\sum_{k= 0,k\neq1}^{\infty} \log \left|\frac{b_{1}-b_{k}}{b_{1}-\bar{b}_{k}}\right|&=&\log\left(\dfrac{j_{\nu+1,1}^2}{2\delta/\kappa_\al^2+j_{\nu+1,1}^2}\right)+\sum_{k=2}^{\infty} \log \left(\frac{\left( j_{\nu+1, k}\right)^{2}-\left( j_{\nu+1, 1}\right)^{2}}{2 \delta / \kappa_\alpha^{2}+\left( j_{\nu+1, 1}\right)^{2}+\left(j_{\nu+1, k}\right)^{2}}\right)\notag\\
		&\leq& \log\left(\dfrac{j_{\nu+1,1}^2}{2\delta/\kappa_\al^2+j_{\nu+1,1}^2}\right)+ \sum_{k=2}^{\infty} \frac{1}{j_{\nu+1, k+1}-j_{\nu+1, k}} \int_{j_{\nu+1, k}}^{j_{\nu+1, k+1}} \log \left(\frac{ x^{2}}{2 \delta / \kappa_\alpha^{2}+ x^{2}}\right) \mathrm{d} x \notag \\ 
		&\leq& \log\left(\dfrac{j_{\nu+1,1}^2}{2\delta/\kappa_\al^2+j_{\nu+1,1}^2}\right)+\frac{1}{\pi} \int_{j_{\nu+1, 2}}^{\infty} \log \left(\frac{ x^{2}}{2 \delta / \kappa_\alpha^{2}+x^{2}}\right) \mathrm{d} x,\label{apoyo6}\\
		&=& \log\left(\dfrac{j_{\nu+1,1}^2}{2\delta/\kappa_\al^2+j_{\nu+1,1}^2}\right) -\frac{j_{\nu+1, 2}}{\pi} \log \left(\frac{1}{1+2 \delta /\left(\kappa_\alpha j_{\nu+1, 2}\right)^{2}}\right)- \frac{2\sqrt{2 \delta}}{\pi\kappa_\alpha}\tan ^{-1}\left(\frac{\sqrt{2 \delta}}{\kappa_\alpha j_{\nu+1, 2}}  \right) ,\notag
	\end{eqnarray}
	where we have used Lemma \ref{consec} and made the change of variables
	$$
	\tau=\frac{ \kappa_\alpha}{\sqrt{2 \delta}} x.
	$$
	%On the other hand, straightforward computations using (73) give for some $C>0$\frac{\sqrt{2 \delta}}{\kappa_\alpha}- 
	From (\ref{logF4}) we get the estimate
	\begin{equation}\label{apoyo7}
		\frac{\Im\left(b_{1}\right)}{\pi} \int_{-\infty}^{\infty} \frac{\log |F(s)|}{\left|s-b_{1}\right|^{2}} \mathrm{~d} s \leq \frac{\delta T }{2}+\log\pts{(\mathcal{K}+\varepsilon) T^{1 / 2}\frac{\abs{J_\nu\pts{j_{\nu+1,1}}}}{\pts{2\kappa_\al}^{1/2}}}.
	\end{equation}
	From (\ref{ena4}), (\ref{aprep4}), (\ref{apoyo6}), and (\ref{apoyo7}), we have
	\begin{eqnarray*}
		\frac{2\sqrt{2 \delta}}{\pi\kappa_\alpha}\tan ^{-1}\left(\frac{\sqrt{2 \delta}}{\kappa_\alpha j_{\nu+1, 2}}\right) -\frac{j_{\nu+1, 2}}{\pi} \log \left(1+\frac{2 \delta}{ \left(\kappa_\alpha j_{\nu+1, 2}\right)^{2}}\right) 
		-(\lambda_1+\delta) T &\leq& \log(\mathcal{K}+\varepsilon)+\log\pts{\frac{(2\kappa_\al T)^{1/2}(j_{\nu+1,1})^{\nu}}{2^{\nu}\Gamma(\nu+1)|J_{\nu}(j_{\nu+1,1})|}}\\
		&&+\log\left(\dfrac{j_{\nu+1,1}^2}{2\delta/\kappa_\al^2+j_{\nu+1,1}^2}\right).
	\end{eqnarray*}
	The result follows by taking $\delta=\kappa_\alpha^2 \pts{j_{\nu+1,2}}^2/2$ and then letting $\varepsilon\rightarrow 0^+$.
	
	\section{Control at the right endpoint}\label{rightend}
	\subsection{Upper estimate of the cost of the null controllability}
	Now we show the null controllability of the system (\ref{problem1Right}). Let $u_{0}\in L^{2}_\beta(0,1)$ given as in (\ref{uoSerie}). We set
	\begin{equation}\label{fserie2}
		f(t):=-\sum_{k=0}^{\infty}\frac{b_{k} \mathrm{e}^{-\lambda_{k} T}}{\Phi_{k}(1)} \psi_{k}(t).
	\end{equation}
	
	Since the sequence $\{\psi_k\}_{k\geq 0}$ is biorthogonal to $\{\mathrm{e}^{-\lambda_k(T-t)}\}_{k\geq 0}$, we have
	\begin{equation}\label{control2}
		\Phi_{k}(1)\int_{0}^{T} f(t)\mathrm{e}^{-\lambda_{k}(T-t)} \mathrm{d} t = -b_{k} \mathrm{e}^{-\lambda_{k} T} = -\left\langle u_{0}, \mathrm{e}^{-\lambda_{k} T}\Phi_{k}\right\rangle_\beta,\quad k\geq 0.
	\end{equation}
	
	Let $u\in C\pts{[0,T];L^{2}_{\beta}(0,1)}$ be the weak solution of system (\ref{problem1Right}).  In particular, for $\tau=T$ we take $z^T=\Phi_k$, $k\geq 0$, then (\ref{weaksolP1R}) and (\ref{control2}) imply that $\left\langle u(\cdot, T), \Phi_{k}\right\rangle_{\beta}=0$ for all $k \geq 0$, therefore $u(\cdot, T) \equiv 0$.\\
	
	Finally, we estimate the norm of the control $f$. From (\ref{Phik_1}), (\ref{psiL1}), (\ref{deriLambda}) and (\ref{fserie2}) we get
	\begin{equation*}\label{finty2}
		C(T, \alpha,\delta)^{-1}\|f\|_{\infty} \leq  \frac{|b_0|}{|\Phi_{0}(1)|} + \frac{1}{\sqrt{2\kappa_\al}2^{\nu}\Gamma(\nu+2)}\sum_{k=1}^{\infty} \frac{|j_{\nu+1,k}|^{\nu}}{|J_{\nu}(j_{\nu+1,k})|}\dfrac{|b_{k}|}{\lambda_{k}}\exp\pts{-\frac{T\lambda_k}{2}-\frac{\omega\lambda_k}{2\sqrt{\theta+1}}}.
	\end{equation*}
	
	By using that $\mathrm{e}^{-x}\leq \mathrm{e}^{-r}r^rx^{-r}$ for all $x,r>0,$ the Cauchy-Schwarz inequality, Lemma \ref{reduce} and the fact that $j_{\nu,k}\geq (k-1/4)\pi$ (by (\ref{below})), (\ref{uoSerie}) and $\lambda_1 \leq \lambda_k$, $k\geq 1$, we obtain that
	\begin{eqnarray*}
		C(T, \alpha,\delta)^{-1}\|f\|_{\infty} &\leq& \frac{|b_0|}{|\Phi_{0}(1)|} + \frac{c\kappa_\al^{-\nu-1}}{\Gamma(\nu+2)}\pts{\dfrac{2\nu+1}{4T}}^{(2\nu+1)/4} \exp\pts{-\frac{2\nu+1}{4}}\exp\pts{-\frac{\omega\lambda_1}{2\sqrt{\theta+1}}-\frac{T\lambda_1}{4}}\sum_{k=1}^{\infty}\frac{|b_k|}{\lambda_k}\\
		&\leq& \frac{|b_0|}{|\Phi_{0}(1)|} + \frac{c\kappa_\al^{-\nu-1}}{\Gamma(\nu+2)}\pts{\dfrac{2\nu+1}{4T\mathrm{e}}}^{(2\nu+1)/4} \exp\pts{-\frac{\omega\lambda_1}{2\sqrt{\theta+1}}-\frac{T\lambda_1}{4}}\pts{\sum_{k=1}^{\infty}|b_k|^2}^{1/2},
	\end{eqnarray*}
	and the result follows by (\ref{reduc}).
	
	%%%%%%%
	%\subsection{Lower estimate of the cost of the null controllability}
	
	\subsection{Lower estimate of the cost of the null controllability}
	Once again, we get a lower estimate of the cost $\widetilde{\mathcal{K}}=\widetilde{\mathcal{K}}_{L^2_\beta}(T,\alpha,\beta,\mu)$.
	We set
	\begin{equation}\label{first3} %x^{\pts{1-\al-\beta}/2}J_\nu\pts{j_{\nu,1}x^{\kappa_\al}}
		u_0(x):= \frac{\abs{J_\nu\pts{j_{\nu+1,1}}}}{(2\kappa_\al)^{1/2}}\Fi_1(x),\,x\in(0,1),
		\quad \text{hence}\quad\|u_0\|^2_\beta=\frac{\abs{J_\nu\pts{j_{\nu+1,1}}}^2}{2\kappa_\al}.
	\end{equation}
	For $\varepsilon>0$ small enough, there exists $f\in \widetilde{U}(\al,\beta,\mu,T,u_0)$ such that
	\begin{equation*}\label{inicost3}
		u(\cdot,T)\equiv 0,\quad \text{and}\quad \|f\|_{L^2(0,T)}\leq (\widetilde{\mathcal{K}}+\varepsilon)\|u_0\|_\beta.
	\end{equation*}
	Then, in (\ref{weaksolP1R}) we set $\tau=T$ and  take $z^\tau=\Fi_k$, $k\geq 0$,  to obtain
	\[\mathrm{e}^{-\lambda_k T}\left\langle u_{0},\Fi_k\right\rangle_\beta=\left\langle u_{0}, S(T)\Fi_k\right\rangle_{\beta}= -
	\int_{0}^{T} f(t) \lim_{x\to 1^{-}}S(T-t) \Fi_k \mathrm{d} t= -\mathrm{e}^{-\lambda_k T}\Fi_k(1)\int_{0}^{T} f(t) \mathrm{e}^{\lambda_k t} \mathrm{d} t.\]
	From (\ref{Phik_1}) and (\ref{first3}) it follows that
	\begin{equation}\label{orto3}
		\int_{0}^{T} f(t) \mathrm{e}^{\lambda_k t} \mathrm{d} t= -\frac{\abs{J_{\nu}\pts{j_{\nu+1, 1}}}}{2\kappa_\al}\delta_{1,k},\quad k\geq 0.
	\end{equation}
	
	Consider the entire function $v: \mathbb{C} \rightarrow \mathbb{C}$ given by
	\begin{equation*}\label{vanalitica}
		v(s):=\int_{-T / 2}^{T / 2} f\left(t+\frac{T}{2}\right) \mathrm{e}^{-i s t} \mathrm{~d} t, \quad s \in \mathbb{C} .
	\end{equation*}
	Therefore,
	\begin{equation}\label{uve3}
		|v(s)| \leq \mathrm{e}^{T|\Im(s)|/2}\int_{0}^{T}|f(t)| \mathrm{d} t \leq (\widetilde{\mathcal{K}}+\varepsilon) T^{1/2}\mathrm{e}^{T|\Im(s)|/2} \left\|u_{0}\right\|_{\beta}.
	\end{equation}
	Moreover, (\ref{orto3}) implies that
	\[v(i\lambda_k)=0\quad\text{for all }k\geq 0, k\neq 1,\quad \text{and}\quad v(i\lambda_1)=-\frac{\abs{J_{\nu}\pts{j_{\nu+1, 1}}}}{2\kappa_\al}\mathrm{e}^{-\lambda_1 T/2}.\]
	Consider the entire function $F(z)$ given by
	\begin{equation}\label{entire3}
		F(s):=v\left(s-i \delta\right), \quad s \in \mathbb{C},\quad\text{with  }\delta=\kappa_\alpha^2 \pts{j_{\nu+1,2}}^2/2.
	\end{equation}
	Clearly, 
	\[ F\left(b_{k}\right)=0, \quad k\geq 0, k\neq 1, \quad \text {where} \quad b_{k}:=i\left(\lambda_k+\delta\right),\quad k\geq 0,\quad\text{and}\]
	\begin{equation}\label{ena3}
		F\left(b_{1}\right)= -\frac{\abs{J_{\nu}\pts{j_{\nu+1, 1}}}}{2\kappa_\al}\mathrm{e}^{-\lambda_1 T/2}.
	\end{equation}
	From (\ref{first3}), (\ref{uve3}) and (\ref{entire3}) we obtain
	\begin{equation}\label{logF3}
		\log |F(s)|\leq \frac{T}{2}|\Im(s)-\delta|+\log\pts{(\widetilde{\mathcal{K}}+\varepsilon) T^{1 / 2}\frac{\abs{J_\nu\pts{j_{\nu+1,1}}}}{\pts{2\kappa_\al}^{1/2}}},\quad s\in\mathbb{C}.
	\end{equation}
	
	We apply Theorem \ref{repreTheo} to the function $F(z)$ given in (\ref{entire3}). Then, (\ref{uve3}) implies that $A\leq T/2$, hence
	\begin{equation}\label{aprep3}
		\log \left|F\left(b_{1}\right)\right|\leq\left(\lambda_1+\delta\right)\frac{T}{2}+\sum_{k=0,k\neq 1}^{\infty} \log \left|\frac{b_{1}-b_{k}}{b_{1}-\bar{b}_{k}}\right|+\frac{\Im\left(b_{1}\right)}{\pi} \int_{-\infty}^{\infty} \frac{\log |F(s)|}{\left|s-b_{1}\right|^{2}} \mathrm{~d}s.
	\end{equation}
	From (\ref{logF3}) we get the estimate
	\begin{equation}\label{apoyo4}
		\frac{\Im\left(b_{1}\right)}{\pi} \int_{-\infty}^{\infty} \frac{\log |F(s)|}{\left|s-b_{1}\right|^{2}} \mathrm{~d} s \leq \frac{T \delta}{2}+\log \left((\widetilde{\mathcal{K}}+\varepsilon) T^{1 / 2} \frac{\left|J_{\nu}\left(j_{\nu+1, 1}\right)\right|}{\pts{2 \kappa_\alpha}^{1/2}}\right).
	\end{equation}
	From (\ref{apoyo6}), (\ref{ena3}), (\ref{aprep3}), and (\ref{apoyo4}) we have
	\begin{equation*}\label{abajo}
		\log\pts{1+\frac{j_{\nu+1,2}^2}{j_{\nu+1,1}^2}}+\pts{\frac{1}{2}-\frac{\log 2}{\pi}}j_{\nu+1,2}-\pts{\lambda_1+\frac{\kappa_\al^2 j_{\nu+1, 2}^2}{2}}T \leq \log(\widetilde{\mathcal{K}}+\varepsilon)+\log (2\kappa_{\al}T)^{1/2},
	\end{equation*}
	the result follows by letting $\varepsilon\rightarrow 0^+$.
	
	\appendix
	\section{Bessel functions}
	We introduce the Bessel function of the first kind $J_{\nu}$ as follows
	\begin{equation}\label{bessel}
		J_{\nu}(x)=\sum_{m \geq 0} \frac{(-1)^{m}}{m ! \Gamma(m+\nu+1)}\left(\frac{x}{2}\right)^{2 m+\nu}, \quad x \geq 0,
	\end{equation}
	where $\Gamma(\cdot)$ is the Gamma function. In particular, for $\nu>-1$ and $0<x \leq \sqrt{\nu+1}$, from (\ref{bessel}) we have (see \cite[9.1.7, p. 360]{abram})
	\begin{equation}\label{asincero}
		J_{\nu}(x) \sim \frac{1}{\Gamma(\nu+1)}\left(\frac{x}{2}\right)^{\nu} \quad \text { as } \quad x \rightarrow 0^{+} .
	\end{equation}
	A Bessel function $J_\nu$ of the first kind solves the differential equation
	\begin{equation}\label{Besselode}
		x^2y''+xy'+(x^2-\nu^2)y=0.
	\end{equation}
	Bessel functions of the first kind satisfy the recurrence formulas (see \cite[9.1.27]{abram}):
	\begin{equation}\label{recur}
		x J_{\nu}^{\prime}(x)-\nu J_{\nu}(x)=-x J_{\nu+1}(x),
	\end{equation}
	\begin{equation}\label{recur2}
		x^{1-\nu}\frac{d}{dx}[x^\nu J_\nu(x)]	=x J_{\nu}^{\prime}(x)+\nu J_{\nu}(x)=x J_{\nu-1}(x).
	\end{equation}
	Recall the asymptotic behavior of the Bessel function $J_{\nu}$ for large $x$, see \cite[Lem. 7.2, p. 129]{komo}.
	\begin{lem}\label{asimxinf}
		For any $\nu \in \mathbb{R}$
		$$
		J_{\nu}(x)=\sqrt{\frac{2}{\pi x}}\left\{\cos \left(x-\frac{\nu \pi}{2}-\frac{\pi}{4}\right)+\mathcal{O}\left(\frac{1}{x}\right)\right\} \quad \text { as } \quad x \rightarrow \infty
		$$
	\end{lem}
	For $\nu>-1$, $\ell,\ell'\in \R$, we have (see \cite[p. 101]{Bowman})
	\begin{equation}\label{orto}
		\int_0^1 x J_\nu(\ell x) J_\nu\left(\ell^{\prime} x\right) \mathrm{d}x=\frac{\ell^{\prime} J_\nu(\ell ) J_\nu^{\prime}\left(\ell^{\prime} \right)-\ell J_\nu\left(\ell^{\prime} \right) J_\nu^{\prime}(\ell )}{\ell^2-\ell^{\prime 2}}.
	\end{equation}
	
	For $\nu >-1$ the Bessel function $J_{\nu}$ has an infinite number of real zeros $0<j_{\nu, 1}<j_{\nu, 2}<\ldots$, all of which are simple, with the possible exception of $x=0$. In \cite[Proposition 7.8]{komo} we can find the next information about the location of the zeros of the Bessel functions $J_{\nu}$:
	\begin{lem}\label{consec}Let $\nu \geq 0$.\\
		1. The difference sequence $\left(j_{\nu, k+1}-j_{\nu, k}\right)_{k}$ converges to $\pi$ as $k \rightarrow\infty$.\\
		2. The sequence $\left(j_{\nu, k+1}-j_{\nu, k}\right)_{k}$ is strictly decreasing if $|\nu|>\frac{1}{2}$, strictly increasing if $|\nu|<\frac{1}{2}$, and constant if $|\nu|=\frac{1}{2}$.\\
	\end{lem}
	
	\begin{proposition}\label{basis}
		Let $\nu>-1$, $0\leq\al<2$ and $\beta\in \R$. The family 
		\[\Fi_0(x) = \sqrt{2(\nu+1)\kappa_\alpha}\,\,x^{(1-\al-\beta)/2+\kappa_\alpha \nu}, \quad\Fi_k(x) = \dfrac{\sqrt{2\kappa_\al}}{|J_\nu(j_{\nu+1,k})|}x^{(1-\al-\beta)/2}J_\nu\pts{j_{\nu+1,k}x^{\kappa_\al}}, k\geq 1,\]
		is an orthonormal basis for $L^2_\beta(0,1)$.
	\end{proposition}
	\begin{proof}
		By using (\ref{recur}) and (\ref{orto}) with $\ell'=j_{\nu+1,k}$, we get
		\[\int_0^1 x J_\nu(\ell x) J_\nu\left(j_{\nu+1,k} x\right) \mathrm{d}x=\frac{\ell J_{\nu+1}(\ell ) J_\nu(j_{\nu+1,k})}{(\ell+j_{\nu+1,k})(\ell-j_{\nu+1,k})}.\]
		By taking the limit as $\ell$ goes to $j_{\nu+1,k}$, and by using (\ref{recur2}) (with $\nu+1$ instead of $\nu$), we obtain
		\begin{equation}\label{relort}
			\int_0^1 x |J_\nu\left(j_{\nu+1,k} x\right)|^2 \mathrm{d}x=\frac{1}{2}J'_{\nu+1}(j_{\nu+1,k})J_\nu(j_{\nu+1,k})=\frac{|J_\nu(j_{\nu+1,k})|^2}{2},\quad k\geq 1.
		\end{equation}
		
		Next, we introduce the following family
		\begin{equation}\label{Dini}
			\Theta_0(x):= \sqrt{2(\nu+1)}x^{1/2+\nu},\quad \Theta_k(x):=\dfrac{\sqrt{2}}{|J_\nu(j_{\nu+1,k})|}x^{1/2}J_\nu\pts{j_{\nu+1,k}x},\quad k\geq 1.
		\end{equation}
		
		In \cite{Hoch} was proved that $\{\Theta_k\}_{k\geq 0}$ is a complete system in $L^2(0,1)$.\\
		
		Then, (\ref{recur}), (\ref{orto}) and (\ref{relort}) imply that $\langle \Theta_k, \Theta_\ell\rangle=\delta_{k,\ell}$ for all $k,\ell\geq 1$. On the other hand, from (\ref{recur2}) with $\nu+1$ instead of $\nu$, we obtain that
		\[(j_{\nu+1,k})^{\nu+2}\int_0^1 x^{\nu+1}J_\nu(j_{\nu+1,k}x)\mathrm{d}x=y^{\nu+1}J_{\nu+1}(y)|_{y=0}^{y=j_{\nu+1,k}}=0,\quad k\geq 1.\]	
		Therefore $\langle \Theta_k, \Theta_0\rangle=0$ for all $k\geq 1$. In conclusion, $\{\Theta_k\}_{k\geq 0}$ is an orthonormal basis for $L^2(0,1)$.\\
		
		Let $\U$ be the unitary operator $\mathcal{U}:L^2(0,1)\rightarrow L^2_\beta(0,1)$ given by
		\begin{equation*}%\label{unitary}
			\U u(x):=\kappa_\al^{1/2}x^{-\al/4-\beta/2}u(x^{\kappa_\al}), \quad u\in L^2(0,1).
		\end{equation*}
		
		Notice that $\U\Theta_k=\Fi_k$, $k\geq 0$, therefore $\Fi_k$, $k\geq 0$, is an orthonormal basis for $L^2_\beta(0,1)$.
	\end{proof}
	
	For $\nu \geq 0$ fixed, we consider the next asymptotic expansion of the zeros of the Bessel function $J_{\nu}$, see\cite[Section 15.53]{Watson},
	\begin{equation}\label{asint}
		j_{\nu, k}=\left(k+\frac{\nu}{2}-\frac{1}{4}\right) \pi-\frac{4 \nu^{2}-1}{8\left(k+\frac{\nu}{2}-\frac{1}{4}\right) \pi}+O\left(\frac{1}{k^{3}}\right), \quad \text { as } k \rightarrow\infty
	\end{equation}
	
	In particular, we have
	\begin{equation}\label{below}
		\begin{aligned}
			&j_{\nu, k} \geq\left(k-\frac{1}{4}\right) \pi \quad \text { for } \nu \in\left[0, 1/2\right], \\
			&j_{\nu, k} \geq\left(k-\frac{1}{8}\right) \pi \quad \text { for } \nu \in\left[1/2,\infty\right].
		\end{aligned}
	\end{equation}
	
	\begin{lem}\label{reduce} For any $\nu > -1$ and any $k\geq 1$ we have
		$$
		\sqrt{j_{\nu+1, k}}\left|J_{\nu}\left(j_{\nu+1, k}\right)\right|=\sqrt{\frac{2}{\pi}}+O\left(\frac{1}{j_{\nu+1, k}}\right)\quad \text{as}\quad k \rightarrow \infty.
		$$
	\end{lem}
	The proof of this result follows by using  (\ref{asimxinf}).
	
	\begin{lem} Let $0\leq \al < 2$, $\beta\in\R$, $a$ and $\nu=\nu(\alpha,\beta,\mu)$ given in (\ref{a_const}) and (\ref{Nu}) respectively, then the following limits are finite
		\begin{equation}\label{lim2}
			\d\OO_{a}(\Fi_0)= \sqrt{2-\al+2\sqrt{\mu(\al+\beta)-\mu}},\quad\,
			\d\OO_{a}(\Fi_k)= \frac{(2\kappa_\al)^{1/2}\pts{j_{\nu+1,k}}^{\nu}}{2^{\nu}\Gamma(\nu+1)\abs{J_{\nu}\pts{j_{\nu+1, k}}}}, \quad k\geq1.
		\end{equation}
	\end{lem}
	\begin{proof}
		This result follows from (\ref{bessel}).
	\end{proof}
	
	%We recall the following representation theorem, see \cite[p. 56]{koosis}.
	%%\begin{theorem}\label{repreTheo} Let $g(z)$ be an entire function of exponential type and assume that
	%%$$
	%%\int_{-\infty}^{\infty} \frac{\log ^{+}|g(x)|}{1+x^{2}} \mathrm{d}x<\infty.
	%%$$
	%%Let $\left\{b_{\ell}\right\}_{\ell \geq 1}$ be the set of zeros of $g(z)$ in the upper half plane $\Im(z)>0$ (each zero being repeated as many times as its multiplicity). Then,
	%%$$
	%%\log |g(z)|=A \Im(z)+\sum_{\ell=1}^{\infty} \log \left|\frac{z-b_{\ell}}{z-\bar{b}_{\ell}}\right|+\frac{\Im(z)}{\pi} \int_{-\infty}^{\infty} \frac{\log |g(s)|}{|s-z|^{2}} \mathrm{d}s,\quad\Im(z)>0,
	%%$$
	%%where
	%%$$
	%%A=\limsup _{y \rightarrow\infty} \frac{\log |g(i y)|}{y} .
	%%$$
	%%\end{theorem}

\end{document}